\documentclass[11pt]{article}
\usepackage{amscd, amsmath, amssymb}

\title{The curvature of a Hessian metric}
\author{Burt Totaro}
\date{  }

\def\Z{\text{\bf Z}}

\def\R{\text{\bf R}}
\def\C{\text{\bf C}}
\def\H{\text{\bf H}}

\def\arrow{\rightarrow}

\def\qed{\ QED }

\def\d{\partial}
\def\gl{\mathfrak gl}

\begin{document}
\maketitle

\newtheorem{theorem}{Theorem}[section]
\newtheorem{corollary}[theorem]{Corollary}
\newtheorem{lemma}[theorem]{Lemma}
\newtheorem{conjecture}[theorem]{Conjecture}
\newtheorem{question}[theorem]{Question}

Given a smooth function $f$ on an open subset of a real vector space,
one can define the associated ``Hessian metric''
using the second derivatives of $f$,
$$g_{ij}:=\partial^2f/\partial x_i\partial x_j.$$
In this paper, inspired by P.M.H.~Wilson's paper on sectional
curvatures of K\"ahler moduli \cite{Wilson},
we concentrate on the case where $f$ is a homogeneous
polynomial (also called a ``form'') of degree $d$ at least 2.
Following Okonek and van de Ven \cite{OV}, Wilson considers the ``index cone,''
the open subset where the Hessian matrix of $f$ is Lorentzian
(that is, of signature $(1,*)$) and $f$ is positive.
He restricts the indefinite metric
$-1/d(d-1)\d^2f/\d x_i\d x_j$ to the hypersurface $M:=\{ f=1\}$
in the index cone, where it is a Riemannian metric, which he calls
the Hodge metric. (In affine differential geometry, this metric
is known as the ``centroaffine metric'' of the hypersurface $M$,
up to a constant factor.)
Wilson considers two main questions about the Riemannian manifold $M$.
First, when does $M$ have nonpositive sectional curvature? (It does have
nonpositive sectional curvature in many examples.) Next, when does $M$
have constant negative curvature?

On the first question, Wilson gave examples of cubic forms $f$ on $\R^3$
to show that the surface $M$ need not have nonpositive curvature
everywhere. But he showed that for every cubic form on $\R^3$ such that $M$
is nonempty (that is, the index cone is nonempty),
$M$ has nonpositive curvature somewhere
(\cite{Wilson}, Prop.\ 5.2). One result of this paper
is to confirm Wilson's suggestion that this statement should fail
for forms of higher degree or on a higher-dimensional space.
Namely, we give examples of a quartic form on $\R^3$ and
a cubic form on $\R^4$ such that $M$ is nonempty and
$M$ has positive sectional curvature on some
2-plane at every point (Lemmas \ref{quartic} and \ref{cubic}).
If Wilson's conjecture that the K\"ahler moduli space of
a K\"ahler manifold has nonpositive sectional curvature is correct,
then these forms cannot occur
as the intersection form on $H^{1,1}(X,\R)$ for a K\"ahler
4-fold with $h^{1,1}=3$, or a K\"ahler 3-fold with $h^{1,1}=4$
(respectively), although they would be allowed by the Hodge index theorem.

Wilson showed that the Riemannian manifold $M$ has constant negative
sectional curvature
when $f$ is a Fermat form $x_1^d-x_2^d-\cdots -x_n^d$
(\cite{Wilson}, Introduction, Example 2).
More generally, we show that $M$ has constant negative curvature $-d^2/4$
when $f$ is a sum of forms of degree $d$ in at most two variables,
$f=\alpha_1(x_1,x_2)+\alpha_2(x_3,x_4)+\cdots$.
The problem of finding forms $f$ such
that the surface $M$ has constant curvature is a special case
of the WDVV equations of string theory, as explained in section
\ref{hess}. In fact, section \ref{hess} lists
a whole series of natural problems
of differential geometry that are essentially equivalent
to the WDVV equations.

The problem of finding all forms $f$ on $\R^3$ such that
the surface $M$ has constant curvature $-d^2/4$ also
has a close relation to classical invariant theory,
in particular to the ``Clebsch covariant'' $S(f)$ studied
by Clebsch \cite{Clebsch} and Dolgachev-Kanev \cite{DK};
see Question \ref{clebsch}.
Using these ideas, we prove that any ternary form $f$ of degree
at most 4 such that
the surface $M$ has constant curvature $-d^2/4$ is in the closure
of the set of forms which can be
written $f=\alpha(x,y)+\beta(z)$ in some linear coordinate system
(Theorem \ref{four}); 
this was known from Wilson's results when $f$ has degree 3
(\cite{Wilson}, Examples, section 5).
Also, we prove a weaker result for
plane curves of any degree: the ternary forms which can be written
$f=\alpha(x,y)+\beta(z)$ in some linear coordinate system always form
an irreducible component of the set of all forms such that
the associated surface $M$ has constant curvature (Theorem \ref{irred}).

Dubrovin showed that Maschke's ternary sextic
$f=x^6+y^6+z^6-10(x^3y^3+y^3z^3+z^3x^3)$, the invariant
of lowest degree for
a complex reflection group of order 648 on $\C^3$,
gives a surface $M$
with constant curvature $-d^2/4$ (\cite{Dubrovinalmost},
Corollary 5.9 and Example 3). 
Dubrovin used the equivalent
statement that this form gives a flat Hessian metric; see
Corollary \ref{flat} below.
The Maschke sextic is not in the closure
of the set of forms which can be written as $\alpha(x,y)+\beta(z)$
in some linear coordinates.
We can ask if the Maschke sextic
is the only ternary form of any degree which gives a surface $M$
of constant curvature $-d^2/4$ while not being in the closure
of the forms $\alpha(x,y)+\beta(z)$ (Question \ref{conjflat}).
It may seem implausible that the Maschke sextic
should be the sole exception here; in some sense, that
would mean that complex reflection groups play a very special role
in this problem. A rough analogy which encourages this belief
is Hertling's theorem: any massive Frobenius manifold
whose Euler vector field has positive degrees arises from
some Coxeter group (\cite{Hertling}, Theorem 5.25).

The paper starts with general formulas for the curvature of Hessian
metrics. We relate Wilson's constructions
to the literature on Hessian metrics, using the notion of
warped products (Lemma \ref{warp}).

I am grateful to Pelham Wilson for many useful conversations,
and to Boris Dubrovin for pointing out the example above, disproving
my original conjecture.

\section{Conventions}
\label{conventions}

Here is Wilson's definition of the ``Hodge metric.'' Let $X$ be a compact
K\"ahler manifold of dimension $d$ at least 2. The cup product on $H^2(X,\Z)$
determines a degree $d$ form $f(\omega):=\omega^d\in \R$
on $H^{1,1}(X,\R)$, and the {\it positive
cone }is defined to be the set of elements $\omega$ of $H^{1,1}(X,\R)$
such that $\omega^d$ is positive. The cup product also determines
an {\it index cone}, as defined by Okonek and van de Ven \cite{OV}:
the set of elements $\omega $ in the positive cone such that the quadratic
form on $H^{1,1}(X,\R)$ defined by $L\mapsto \omega^{d-2}L^2$
has signature $(1,h^{1,1}-1)$. The Hodge index theorem says that
any K\"ahler class in $H^{1,1}(X,\R)$ lies in the index cone.
Let $W_1$ be the set of points $\omega$ in the index cone with
$\omega^d=1$. Then $W_1$ is a smooth manifold, whose tangent space at
a point $\omega$ is the set of $L$ in $H^{1,1}(X,\R)$ such that
$\omega^{d-1}L=0$. The {\it Hodge metric }is the Riemannian metric
on $W_1$ defined by, for tangent vectors $L_1$ and $L_2$ at
a point $\omega $ in $W_1$,
$$(L_1,L_2)=-\omega^{d-2}L_1L_2.$$
One computes easily that this metric is the
restriction to $W_1$ of the pseudo-Riemannian Hessian metric
$-1/d(d-1)\d^2 f/\d x_i\d x_j$ on $H^{1,1}(X,\R)$. In this paper,
we use $-1/d(d-1)\d^2 f/\d x_i \d x_j$ as our basic pseudo-Riemannian
metric on an open subset of a real vector space. We also study
its restriction to the
hypersurface $M:=W_1=\{ f=1 \}$. We usually call this metric
on $\R^n$ the Hessian metric associated
to $f$. Outside this paper, that name is usually restricted
to the metric $\d^2 f/\d x_i \d x_j$ with no constant factor.

We sometimes use the notation $f_i$ for $\d f/\d x_i$,
$f_{ij}$ for $\d^2f /\d x_i\d x_j$, and so on. Starting from section
\ref{const}, we consider ternary forms $f(x,y,z)$, and we
identify these variables with $x_1,x_2,x_3$, so that
$f_{23}$ denotes $\d^2 f/\d y \d z$.

\section{The Hessian metric of a homogeneous polynomial}
\label{hess}

We begin by recalling the formula for the curvature of the Hessian
metric $g_{ij}:=\partial^2f/\partial x_i\partial x_j$
associated to a smooth function $f$ on a domain
in $\R^n$. When $f$
is a homogeneous polynomial, the metric on a domain in $\R^n$
is determined in a simple way, as a warped product, from its
restriction to the hypersurface $M:=\{ f=1\}$. We deduce a formula
for the curvature of $M$, using the known result on $\R^n$;
this approach is simpler than trying to compute the curvature
directly on $M$. 
Finally, we recall the known relation of the metric
on $M$ with another Hessian metric on $\R^n$,
$\d ^2(\log f)/\d x_i\d x_j$, and with the ``centroaffine metric'' of $M$.
Along the way, we relate these constructions to the WDVV equations,
although this is not needed for the rest of the paper.

Hessian metrics have also been called affine K\"ahler metrics
(\cite{CYflat}, \cite{Loftin}, \cite{KS}),
since any K\"ahler metric on a complex manifold has an analogous
local description as $\partial^2f/\partial z_i\partial \overline{z_j}$.
But the name ``Hessian'' seems preferable on historical grounds.

Hessian metrics are
a very natural way to construct Riemannian or pseudo-Riemannian
metrics. For example, they have been
used to define a canonical Riemannian metric on
an arbitrary convex cone \cite{Vinberg}. A variant of Hessian metrics
can be used to define
a canonical Riemannian metric on any convex domain, using the solution
of a Monge-Amp\`ere equation by Loewner-Nirenberg and Cheng-Yau
(\cite{LN}, \cite{CYreg}, \cite{Sasaki}).
There are recent surveys on Hessian metrics by
Duistermaat \cite{DuistermaatHess} and Shima-Yagi \cite{SY}.
Duistermaat observed that one canonical Hessian metric on a convex
domain can have positive curvature (\cite{Duistermaatbound}, 8.4); that
phenomenon is roughly analogous to some examples below
(Lemmas \ref{quartic} and \ref{cubic}).

There is a geometric interpretation of Hessian metrics which has appeared
recently in mirror symmetry (Hitchin \cite{Hitchin},
Leung \cite{Leung}, Kontsevich-Soibelman \cite{KS}). It is worth
mentioning, although it will only be used in this section of the paper.
Given a real vector space $V$,
the cotangent bundle $T^*V=V\oplus V^*$ has a natural symplectic
form $\sum_i dx_i\wedge dp_i$, as is well known, but also 
a natural pseudo-Riemannian metric, $\sum_i dx_idp_i$. Consider
any Lagrangian submanifold of $T^*V$.
Such a submanifold, if its tangent space is in general position,
can be locally viewed as the graph of an exact 1-form $df$, for a
smooth function $f$ on an open subset of $V$. Then the restriction
of the pseudo-Riemannian metric on $T^*V$ to this submanifold
is precisely the Hessian metric associated to $f$. Generically,
a Lagrangian submanifold of $T^*V$ can also be viewed as the graph of
$d\hat{f}$ for a function $\hat{f}$ on an open subset of $V^*$,
known as the Legendre
transform of $f$. So this picture ``explains'' the classical fact
that a function 
and its Legendre transform determine isometric Hessian metrics.

To begin our computations, let $f$ be a smooth function on a region in $\R^n$,
and define $g_{ij}=\partial^2f/\partial x_i\partial x_j$.
We assume that $\det(g_{ij})$ is not zero; then $g_{ij}$ is
a pseudo-Riemannian metric. (O'Neill's book \cite{O'Neill}
is a convenient reference for pseudo-Riemannian metrics.)
The formula for the curvature tensor of a Hessian metric is well known.
The easiest way to compute it uses
the following classical formulas for any pseudo-Riemannian metric,
giving the Christoffel symbols of the
first kind and the
curvature tensor (\cite{Spivak}, Chapter 4.D, equation (***)):
\begin{align*}
\Gamma_{ijk}&= \frac{1}{2}(g_{ij,k}+g_{jk,i}-g_{ik,j})\\
R_{ijkl}&= -\frac{1}{2}(g_{ik,jl}+g_{jl,ik}-g_{il,jk}-g_{jk,il})
-\sum_{p,q} g^{pq}(\Gamma_{jpl}\Gamma_{iqk}-\Gamma_{ipl}\Gamma_{jqk}),
\end{align*}
where $g^{ij}$ is the inverse of the matrix $g_{ij}$, $g_{ij,k}$ means
$\d g_{ij}/\d x_k$, and so on.
For a Hessian metric, it follows immediately that $\Gamma_{ijk}=f_{ijk}/2$
and
$$R_{ijkl}=-\frac{1}{4}\sum_{p,q} g^{pq}(f_{jlp}f_{ikq}-f_{ilp}f_{jkq}).$$
It is remarkable that
the curvature of a Hessian metric depends only on the derivatives
of $f$ to order at most three, whereas one would expect fourth derivatives
of $f$ to come in; Duistermaat gives some explanation for this phenomenon
\cite{DuistermaatHess}.

With our conventions, the sectional curvature
of the 2-plane spanned by $\d/\d x_1$ and $\d/\d x_2$ is
$R_{1212}/(g_{11}g_{22}-g_{12}^2)$.

The second half of
this paper (see Question \ref{conjflat})
will be devoted to the study of flat Hessian metrics;
that is, functions $f$ such that the expression $R_{ijkl}$ is
identically zero. This system of partial differential equations
looks much like the WDVV equations of string theory,
$$\sum_{p,q} (f_{jlp}f_{ikq}-f_{ilp}f_{jkq})=0.$$
In fact, Kito showed that
the problem of finding flat Hessian metrics is precisely equivalent
(by changing to coordinates adapted to the metric) to the WDVV equations
(\cite{Kito}, Lemma 2.2). By the interpretation of Hessian metrics
in terms of Lagrangian submanifolds discussed earlier, it follows
that the WDVV equations also describe the flat Lagrangian submanifolds
of $\R^{2n}$, where $\R^{2n}$ is given
a natural symplectic structure and pseudo-Riemannian metric of
signature $(n,n)$.

In fact, there is a whole series of natural problems of differential
geometry that are essentially equivalent to the WDVV equations.
The problem of classifying flat Lagrangian submanifolds in $\R^{2n}$
with the above pseudo-Riemannian metric is equivalent to the WDVV equations,
as we have just mentioned, but another real form of the same
problem is
the classification of flat Lagrangian
submanifolds of the K\"ahler manifold $\C^n$. Terng showed that the latter
problem is essentially equivalent to a natural integrable system 
of first-order PDEs which
she defined, the $U(n)/O(n)$-system (\cite{Terng}, Prop.\ 3.5.3).
Yet another equivalent problem is that of finding Egorov metrics,
that is, flat metrics of the form: $g_{ij}=0$ for $i\neq j$ and
$g_{ii}=\d f/\d x_i$ (\cite{Terng}, Theorem 3.4.3). Finally,
Ferapontov found that the problem in affine differential geometry
of classifying hypersurfaces in $\R^n$ whose centroaffine metric
is flat also reduces to the WDVV equations \cite{Ferapontov}.
There is a similar reduction
of hypersurfaces whose centroaffine metric has nonzero
constant curvature to WDVV, as follows from Corollary \ref{flat} and the
comments after it.

One basic result about this family of problems is
Moore and Morvan's theorem
that flat Lagrangian submanifolds of $\C^n$ are classified locally by
$n(n+1)/2$ functions of one variable \cite{MM}. It follows, for example,
that flat Hessian metrics
on $\R^n$ are likewise classified locally by $n(n+1)/2$
functions of one variable.

The solutions of the WDVV equations which satisfy certain normalization
and homogeneity
conditions are described geometrically by Dubrovin's theory of
Frobenius manifolds \cite{Dubrovin}. The most important
Frobenius manifolds, the semisimple ones, are locally classified
by finitely many numbers (rather than functions), although non-semisimple
Frobenius manifolds of dimension at least 4
can depend on arbitrary functions of one variable
(\cite{DubrovinPainleve}, Exercise 3.1). Dubrovin also discovered
the relation of Frobenius manifolds with some of the
integrable systems mentioned above, such as Egorov metrics.

Now suppose that the domain $U$ in $\R^n$ is a cone (that is,
$U$ is preserved under  multiplication by positive real
numbers). Suppose also that
$f$ is homogeneous of some degree $d>1$, in the sense that
$$f(\lambda x)=\lambda^d f(x)$$
for all $\lambda >0$. Finally, assume that $f>0$ on $U$.

From here on in the paper, we will use the metric
$g_{ij}=-1/d(d-1)\d^2f/\d x_i\d x_j$ on $U$, to conform
to Wilson's conventions as explained in section \ref{conventions}.
The curvature tensor of
that metric follows from the formula above:
$$R_{ijkl}=-\frac{1}{4d^2(d-1)^2}
 \sum_{p,q} g^{pq}(f_{jlp}f_{ikq}-f_{ilp}f_{jkq}),$$
where $g^{ij}$ now denotes the inverse of the new matrix $g_{ij}$.

Let $M$ be the hypersurface $\{ f=1\}$. Our assumption that
$f$ is homogeneous of degree $>1$ implies that the restriction
of the Hessian metric to
$M$ is nondegenerate, by the calculations in the following lemma.
This lemma shows that
the metric on $M$ determines the metric on
the whole open set $U$ in a simple way. For any pseudo-Riemannian
manifold $M$, define the {\it warped product }$(-1)\R^{>0}\times_t M$ to be
the product manifold $\R^{>0}\times M$ with metric such that
$T(\R^{>0})$ is orthogonal to $TM$ everywhere, the inner product on $TM$ at
$(t,x)\in \R^{>0}\times M$ is $t^2$ times the given inner product, and
the metric on $(-1)\R^{>0}$ is the negative definite
metric $\langle \partial/\partial t,
\partial/\partial t\rangle =-1$. (This construction is most familiar
when $M$
is the sphere $S^{n-1}$ of radius 1: then the warped
product $\R^{>0}\times_t S^{n-1}$ is isometric to
flat Euclidean space $\R^n-0$.)
Also, for any pseudo-Riemannian
manifold $M$, let $aM$ denote the same manifold with inner product
multiplied by $a$, $\langle x,y\rangle_{aM}=a\langle x,y\rangle_M$.

\begin{lemma}
\label{warp}
Give $U$ the pseudo-Riemannian metric $-1/d(d-1)\d ^2 f/\d x_i \d x_j$, and
give the hypersurface $M=\{ f=1\}$ the restricted metric. Then
the map
$$\varphi(t,x):= t^{2/d}x$$
gives an isometry from the warped product $(4/d^2)((-1)\R^{>0}\times_t
 (d^2/4)M)$
to $U$.
\end{lemma}

{\bf Proof. }We use the Euler identity repeatedly: for a homogeneous
function $f$ of degree $d$, we have
$$\sum x_if_i=df,$$
where $f_i$ denote the derivatives of $f$ at a point $x=(x_1,\ldots,x_n)$.
Clearly $\varphi$ gives a diffeomorphism from $\R^{>0}\times M$ to $U$.

Let $v$ be a tangent vector to the hypersurface
$M$ in $U$ at a point $x$. 
Clearly the metric $-1/d(d-1)\partial^2f/\partial x_i\partial x_j$
on $U$ is homogeneous of degree $d-2$. That is,
$$\langle v,v \rangle_{\lambda x}=\lambda^{d-2}
\langle v,v \rangle_{x}$$
for any $\lambda>0$,
where the subscripts denote the point in $U$ where the inner product
is taken. Therefore,
$$\langle \lambda v,\lambda v \rangle_{\lambda x}=\lambda^d
\langle v,v \rangle_{x}.$$
Pulling the metric on $U$ back to $\R^{>0}\times M$ by the map
$\varphi(t,x)=t^{2/d}x$ gives a metric such that
$$ \langle v,v \rangle_{(t,x)}=t^2
\langle v,v \rangle_{(1,x)}.$$
The warped product metric on $(4/d^2)((-1)\R^{>0}\times_t (d^2/4)M)$
has the same homogeneity property, and agrees with the given metric
on $M$ when $t=1$. So the two metric are the same on tangent vectors
in $TM$ at any point $(t,x)$ in $\R^{>0}\times M$.

To check that $T\R$ is orthogonal to $TM$ in the metric pulled back
from $U$, we compute the inner product at a point $x\in U$ of
the outward vector $\sum x_i\partial/\partial x_i$ with
a tangent vector $\sum v_i\partial/\partial x_i$ to the hypersurface
$f=c$, which means that $\sum v_if_i=0$:
\begin{align*}
\langle \sum x_i\partial/\partial x_i, \sum v_i \partial/\partial x_i\rangle
&= -1/d(d-1)\sum_{i,j} x_iv_jf_{ij}\\
&= -1/d\sum_j v_jf_j\\
&= 0.
\end{align*}

Finally, we compute the inner product of the outward vector
$\sum x_i\partial/\partial x_i$ at a point $x\in U$ with itself:
\begin{align*}
\langle \sum x_i\partial/\partial x_i, \sum x_i \partial/\partial x_i\rangle
&= -1/d(d-1)\sum_{i,j} x_ix_jf_{ij}\\
&= -1/d\sum_j x_jf_j\\
&= -f(x).
\end{align*}
We compute that 
$$\frac{\d }{\d t}\varphi(t,y)=\frac{2}{td}\varphi(t,y).$$
Using that, let us compute the length squared of
the tangent vector $\partial/\partial t$ at the point
$(t,y)$ in $\R^{>0}\times M$
in the metric pulled back from $U$. Let
$x=\varphi(t,y)$ in $U$.
\begin{align*}
\langle \partial/\partial t, \partial/\partial t \rangle
&= \frac{4}{d^2t^2}\langle \sum x_i\partial/\partial x_i,
\sum x_i\partial/\partial x_i\rangle \\
&= -\frac{4}{d^2t^2}f(x) \\
&= -\frac{4}{d^2}.
\end{align*}
This agrees with the inner product $\langle
\partial/\partial t, \partial/\partial t \rangle$ in
the warped product metric $(4/d^2)((-1)\R^{>0}\times_t (d^2/4)M)$,
as we want. \qed

Using O'Neill's formula for the sectional curvature of a warped
product of pseudo-Riemannian manifolds, we now deduce a simple
relation between the curvature of the open set $U$ in $\R^n$
and the hypersurface $M$.

\begin{corollary}
\label{oneill}
Let $x$ be a point in the hypersurface $M$ in $U$, $P$
a nondegenerate 2-plane in the tangent space to $M$ at $x$, and $u$
a positive real number. Let $K_M(P)$ be the sectional curvature
of $M$ at the 2-plane $P$. Then the sectional curvature of $U$
at the point $cx$ and the 2-plane $cP$ is
$$K_U(cP)=\frac{1}{c^d}(K_M(P)+d^2/4).$$
\end{corollary}

{\bf Proof. }This follows from Lemma \ref{warp} together with O'Neill's
curvature formula for a warped product (\cite{O'Neill}, Proposition 7.42).
\qed

\begin{corollary}
\label{flat}
The pseudo-Riemannian Hessian metric $-1/d(d-1)\partial^2f/\partial
x_i\partial x_j$ on $U$ is flat if and only
if its restriction to the hypersurface $M:=\{ f=1\}$
has constant sectional curvature $-d^2/4$.
\end{corollary}

{\bf Proof. }Suppose $M$ has constant sectional curvature $-d^2/4$.
Then $-(d^2/4)M$ has constant sectional curvature $1$, and so
the warped product $\R^{>0}\times_t (-d^2/4)M$ is flat by O'Neill's formulas
(\cite{O'Neill}, Proposition 7.42).
(The formulas are the same as those showing that the warped product
$\R^{>0}\times_t S^{n-1}$ is isometric to $\R^n-0$.) By Lemma
\ref{warp}, it follows that $U$ is flat. The converse is immediate
from Corollary \ref{oneill}. \qed

Wilson already noticed that the value $-d^2/4$ for the
curvature of $M$ plays a special role, and Corollary \ref{flat} provides
an explanation of this phenomenon.

To conclude this section, I will state the relation
between the
Hessian metric $\partial^2(\log f)/\partial x_i\partial x_j$ on
$U$ and the above metric on the hypersurface $M$, which is proved
by the same kind of calculation as Lemma \ref{warp}
(Loftin \cite{Loftin}, Theorem 1). Loftin also mentions that
the above metric on $M$
(with a different normalization, namely
the metric $-1/d\; \d^2f/\d x_i\d x_j$ restricted to $M$)
is known in affine differential geometry
as the {\it centroaffine metric }of $M$. This is not needed for the
rest of the paper, but the Hessian metric
associated to the logarithm of a homogeneous function is used in
many papers on convex cones (\cite{Vinberg}, \cite{Guler96},
\cite{Guler97}).

\begin{lemma}
Give $M=\{ f= 1\}$ the pseudo-Riemannian metric obtained by
restricting the Hessian metric $-\partial^2f/\partial x_i\partial x_j$
on $U$. Then the map
$$\alpha(t,x):=e^{t/\sqrt{d}}x$$
is an isometry from the product $\R\times M$ to the
Hessian metric $-\partial^2(\log f)\partial x_i\partial x_j$ on $U$.
Here the real line $\R$ is given
its usual Riemannian metric.
\end{lemma}

{\bf Example. }The most famous examples of Hessian metrics
are those associated to the ``symmetric cones''. Namely, let $f$ be
either a Lorentzian quadratic form $x_1^2-x_2^2-\cdots -x_n^2$,
or else the determinant function on the real vector space of
$n\times n$ symmetric, Hermitian, quaternion Hermitian,
or octonion Hermitian matrices; in the octonion case, we set $n=3$.
Restrict $f$ to the convex cone $\{ f>0,x_1>0\} $ in
the quadratic form case, and to the cone of positive definite matrices
in the other cases. Then the Hessian metric $-1/d(d-1)\d^2f/\d x_i\d x_j$
restricted to the hypersurface $M:=\{ f=1\}$ in this convex cone
gives a Riemannian symmetric space of noncompact type.
In the quadratic form case,
$M$ is hyperbolic space of dimension $n-1$; in the other cases,
we get the symmetric spaces $SL(n,\R)/SO(n)$, $SL(n,\C)/SU(n)$,
$SL(n,\H)/Sp(n)$, and $E_6/F_4$. A reference on symmetric cones
and more general homogeneous convex cones is Vinberg \cite{Vinberg}.

\section{Hessian metrics on $\R^3$ and the Clebsch covariant}

Wilson gave a simple formula for the curvature of the Hessian
metric associated to a cubic form $f$ on $\R^3$ (\cite{Wilson},
Theorem 5.1). (More precisely, he considered the curvature of the
surface $\{ f=1 \}$ in $\R^3$, but that is equivalent in
view of Corollary \ref{oneill}.)
In this section, we extend his formula to apply to
a homogeneous function $f$ of any degree greater than 2
on $\R^3$. The key ingredient in Wilson's formula is the
classical Aronhold invariant $S$ of a plane cubic curve;
our formula involves a natural generalization, the Clebsch covariant.

For a cubic in Weierstrass form, the Aronhold invariant is $S=2^23^{-3}g_2$;
that is, $S$ is a constant multiple of the Eisenstein series of weight 4.
Explicitly, for any cubic form
$$f=a_3x_3^3+3(a_2x_1+b_2x_2)x_3^2+3(a_1x_1^2+2b_1x_1x_2+c_1x_2^2)x_3
+(a_0x_1^3+3b_0x_1^2x_2
+3c_0x_1x_2^2+d_0x_2^3),$$
the Aronhold invariant is given by
\begin{multline*}
S=-(a_0a_2-a_1^2)c_1^2+(a_0a_3-a_1a_2)c_0c_1-(a_1a_3-a_2^2)c_0^2
-b_0^2a_3c_1+b_0b_1(3a_2c_1+a_3c_0) \\
-(b_0b_2+2b_1^2)(a_1c_1+a_2c_0)+b_1b_2(a_0c_1+3a_1c_0)-b_2^2a_0c_0\\
+d_0[b_0(a_1a_3-a_2^2)-b_1(a_0a_3-a_1a_2)+b_2(a_0a_2-a_1^2)]
+(b_0b_2-b_1^2)^2.
\end{multline*}
Here I changed the sign of the Aronhold invariant from Elliott
(\cite{Elliott}, p.~377) to agree with the convention in Wilson's paper
\cite{Wilson}.

Clebsch observed that any invariant for forms
of a given degree (in a given number of variables) extends in a natural
way to a covariant for forms of any larger degree (in the same
number of variables) \cite{Clebsch}.
One way to describe the construction is that
we view the given invariant as an $SL(n)$-equivariant differential operator.
For example, this procedure turns the discriminant invariant of
quadratic forms into the Hessian covariant for forms of any degree.

This procedure turns the Aronhold invariant of cubic forms into the
{\it Clebsch covariant }$S(f)$ for forms $f(x_1,x_2,x_3)$ of any degree,
defined by Clebsch \cite{Clebsch} and further studied by Dolgachev and Kanev
\cite{DK}.
Explicitly,
$S(f)$ is defined by the same formula as the Aronhold invariant
(above), but with $a_i,b_i,c_i,d_i$ defined as:
$$\begin{array}{llll}
a_3=f_{333}\\
a_2=f_{133}& b_2=f_{233}\\
a_1=f_{113}& b_1=f_{123} & c_1=f_{223}\\
a_0=f_{111} & b_0=f_{112} & c_0=f_{122} & d_0=f_{222}.
\end{array}$$
Here $f_{ijk}$ denotes the third derivative of $f$ with respect to
the variables $x_i,x_j,x_k$. This definition has the property
that the Clebsch covariant $S(f)$ of a cubic form $f$ is
$2^43^4$ times the Aronhold invariant $S$ of $f$. In general,
the Clebsch covariant of a form $f$ of degree $d$ is a form
of degree $4(d-3)$.

We now give our simplified
formula for the curvature of the Hessian metric associated to
a homogeneous function on $\R^3$. 
(To compare our formula
with Wilson's formula on cubic forms $f$ (\cite{Wilson},
Theorem 5.1), one has to remember the above factor of $2^43^4$.)

\begin{theorem}
\label{curv}
Let $f$ be a smooth homogeneous function
of degree $d>2$ on an open subset $U$ of $\R^3$. Suppose that
the Hessian determinant of $f$ is nonzero on $U$, and
consider the
pseudo-Riemannian Hessian metric
$-1/d(d-1)\partial^2 f/\partial x_i\partial x_j$ on $U$.
Let $M=\{ f=1\}$. Then the sectional curvature of $U$
on the tangent 2-plane to $M$ at a point is given in terms of
the Hessian determinant $H(f)$ and the Clebsch covariant $S(f)$ by
$$K=\frac{d^2(d-1)^2}{4(d-2)^2}\frac{S(f)f}{H(f)^2}.$$
The sectional curvature of the surface $M$ at the same point is
$$K_M=-\frac{d^2}{4}+\frac{d^2(d-1)^2}{4(d-2)^2}\frac{S(f)f^2}{H(f)^2}.$$
\end{theorem}

Notice that the theorem considers a point on $M$, thus with $f=1$,
and so the factors of $f$ in the formulas are not strictly necessary.
We include them so that the first formula gives, more generally, the
sectional curvature of $U$ on the tangent 2-plane to a point in any
level set $\{ f=\lambda\}$. This sectional curvature
is homogeneous of degree $-d$ by Corollary \ref{oneill}. On the other hand,
the formula for $K_M$ includes $f^2$ so as to be homogeneous of degree 0:
this formula gives, at any point $x$ of $U$, the sectional curvature
of the surface $M$ at the point of $M$ which is a scalar multiple of $x$.

{\bf Proof. }By Corollary \ref{oneill}, the second formula follows
from the first. So we just need to compute the sectional curvature
of $U$ at the tangent 2-plane to $M$ at a point.

The equality we want is an algebraic identity among the derivatives
of $f$. So it suffices to prove the same identity for
holomorphic functions $f$ on an open subset $U$ of $\C^3$ which
are homogeneous of degree $d$. In this situation, the ``metric''
$-1/d(d-1)\d^2 f/\d x_i\d x_j$ is a nondegenerate symmetric bilinear
form on the tangent bundle to $U$. Sectional curvature is defined
by the same formulas as for real metrics.

Both the sectional curvature of $U$ at the tangent 2-plane to $M$
and the right side of the above formula are
unchanged under a change of coordinates in $GL(3,\C)$. For the
sectional curvature, this is clear. For the right side, it follows
from the identities, for $A\in GL(3,\C)$:
\begin{align*}
H(fA)|_x &= H(f)|_{Ax}\det(A)^2 \\
S(fA)|_x &= S(f)|_{Ax}\det(A)^4.
\end{align*}

So we can
assume that the given point of $M$ is $(0,0,1)$. Since this point is in $M$,
we have $f(0,0,1)=1$. By the Euler identity, it follows that
$f_3|_{(0,0,1)}=\sum x_if_i|_{(0,0,1)}=d\cdot f(0,0,1)=d$.
By a further change of coordinates in $GL(3,\C)$,
we can assume that the tangent plane to $M$ at $(0,0,1)$ is
spanned by the vectors $\d/\d x_1$ and $\d/\d x_2$. Equivalently,
$\d f/\d x_1=\d f/\d x_2=0$ at the point $(0,0,1)$. By the Euler identity
again, at the point $(0,0,1)$ we have $f_{13}=f_{23}=0$ and $f_{33}=d(d-1)$.
Finally, we are assuming that the Hessian determinant of $f$ is nonzero
at the point $(0,0,1)$.
So we can make one last change of coordinates in
$GL(3,\C)$ so as to make $f_{11}=d(d-1)$, $f_{12}=0$, and
$f_{22}=d(d-1)$ at the point $(0,0,1)$. (Over $\R$ we would have
several cases, depending on the signature of the Hessian metric.)

By the Euler identity, we have $f_{113}=d(d-1)(d-2)$, $f_{123}=0$,
$f_{223}=d(d-1)(d-2)$, $f_{133}=0$, $f_{233}=0$, and
$f_{333}=d(d-1)(d-2)$ at the point $(0,0,1)$. Also, the matrix
$g_{ij}=-1/d(d-1)\d^2 f/\d x_ix_j$ is the matrix $-1$
at the point $(0,0,1)$, and so
its inverse $g^{ij}$ is also the matrix $-1$. We can now
use the formula in section \ref{hess}
for the curvature tensor of the Hessian metric
$-1/d(d-1)\d^2f/\d x_i\d x_j$ on $\C^3$, applied to the 2-plane
$T_{(0,0,1)}M$ spanned by $\d/\d x_1$ and $\d/\d x_2$:
\begin{align*}
R_{1212} &= -\frac{1}{4d^2(d-1)^2}\sum_{p,q}g^{pq}(f_{11q}f_{22p}
-f_{12p}f_{12q})\\
&= \frac{1}{4d^2(d-1)^2}\sum_p (f_{11p}f_{22p}-f_{12p}^2)\\
&= \frac{1}{4d^2(d-1)^2}(f_{111}f_{122}-f_{112}^2+f_{112}f_{222}
-f_{122}^2+d^2(d-1)^2(d-2)^2).
\end{align*}
The sectional curvature of $U$ at this 2-plane is
$R_{1212}/(g_{11}g_{22}-g_{12}^2)$. This is equal to $R_{1212}$, and so
the sectional curvature is given by the formula above.

We compare this to the Clebsch covariant of $f$ at the same
point $(0,0,1)$. Most terms vanish because of what we know
about the third derivatives of $f$ (in the notation used
in defining the Clebsch covariant, $a_2=f_{133}$, $b_2=f_{233}$,
and $b_1=f_{123}$ are zero). What remains is:
\begin{align*}
S(f)&= a_1^2c_1^2+a_0a_3c_0c_1-a_1a_3c_0^2-b_0^2a_3c_1+d_0b_0a_1a_3\\
&= d^2(d-1)^2(d-2)^2(f_{111}f_{122}-f_{112}^2+f_{112}f_{222}
-f_{122}^2+d^2(d-1)^2(d-2)^2).
\end{align*}
Therefore the curvature of $U$ at the above 2-plane is
$S(f)/(4d^4(d-1)^4(d-2)^2)$. At this point $(0,0,1)$, $f$
is equal to 1 and the Hessian determinant of $f$
is $d^3(d-1)^3$. So we can rewrite the
curvature of $U$ at the above 2-plane as
$$K=\frac{d^2(d-1)^2}{4(d-2)^2}\frac{S(f)f}{H(f)^2}.$$
This is the formula we want. \qed

\section{A quartic form on $\R^3$}
\label{sectquartic}

In this section, we answer a question 
raised in the introduction to Wilson's paper,
by giving the first example of a real form $f$ such that the
submanifold $M=\{ f =1\}$ of the index cone is nonempty and has
positive curvature everywhere. Here $M$ is given the Riemannian
metric $-1/d(d-1)\d^2f/\d x_i\d x_j$, and the index cone is defined
in section \ref{conventions}.
The example here is a quartic form
in $\R^3$, while the next section gives a cubic form on $\R^4$ with
analogous properties. These examples are optimal: $M$
has constant curvature $-1$ whenever $f$ is a quadratic form
with nonempty index cone, and $M$ has at least one point of
nonpositive curvature when $f$ is a cubic form on $\R^3$
with nonempty index cone,
by Wilson (\cite{Wilson}, Prop.\ 5.2).

\begin{lemma}
\label{quartic}
For the real quartic form $f=xyz(x+y+z)$, the index cone in $\R^3$
is nonempty, and the surface $M=\{ f=1 \}$ inside the index cone
has positive curvature everywhere.
\end{lemma}

Note that
the form $f$ is not equivalent under $GL(3,\R)$
to its negative, and indeed our argument will not apply
to the negative of $f$.

{\bf Proof. }
The Hessian matrix of $f$ is
$$(\partial^2f/\partial x_i\partial x_j)=\begin{pmatrix}
2yz & 2xz+2yz+z^2 & 2xy+2yz+y^2\\ 2xz+2yz+z^2 & 2xz & 2xy+2xz+x^2 \\
2xy+2yz+y^2 & 2xy+2xz+x^2 & 2xy \end{pmatrix}.$$
We compute that the point $(1,1,1)$ is in the index cone;
that is, $f$ is positive at this point and the Hessian matrix
has signature $(+,-,-)$ at this point. So the index cone is nonempty.
More precisely, we compute that
the Hessian of $f$ (the determinant of the Hessian matrix) is
$$H(f)=6xyz(x+y+z)(x^2+y^2+z^2+xy+xz+yz).$$
Here the quadratic form $x^2+y^2+z^2+xy+xz+yz$ on $\R^3$ is positive
definite. Using that, it is straightforward to check that the index cone of $f$
is equal to the ``positive cone'' $\{ f>0 \}$.

We now show that the surface $M$,
with the metric $-1/d(d-1)\d^2f/\d x_i\d x_j$,
has positive curvature everywhere.
By Theorem \ref{curv}, for a quartic form $f$ in 3 variables,
the curvature of $M$ has the form:
$$K=-4+9 \frac{S(f)f^2}{H(f)^2}.$$
Thus, we need to show that $S(f)>2^23^{-2}H(f)^2/f^2$
on the whole index cone. By the above calculation of the Hessian,
it is equivalent to show that the quartic form
$$S(f)-2^{4}(x^2+y^2+z^2+xy+xz+yz)^2$$
is positive on the whole index cone.

We compute that
\begin{multline*}
S(f)=2^{4}(x^4+2x^3y+3x^2y^2+2xy^3+y^4+2x^3z+7x^2yz+7xy^2z+2y^3z+3x^2z^2\\
+3y^2z^2+7xyz^2+2xz^3+2yz^3+z^4).
\end{multline*}
Therefore the above difference is
\begin{align*}
S(f)-2^{4}(x^2+y^2+z^2+xy+xz+yz)^2 &=2^43xyz(x+y+z)\\
&= 2^43f.
\end{align*}
Since $f$ is positive on the index cone, the above difference is positive on
the index cone, as we want. \qed

\section{A cubic form on $\R^4$}

This section gives another example which answers the question
raised in Wilson's paper and discussed in section \ref{sectquartic}.
The example here is a cubic form on $\R^4$
rather than a quartic form on $\R^3$.

\begin{lemma}
\label{cubic}
For the real cubic form $f=(x_0^2+x_1^2-x_2^2-x_3^2)x_3$,
the index cone in $\R^4$
is nonempty, and the 3-manifold $M=\{ f=1 \}$ inside the index cone
has positive sectional curvature on some 2-plane at every point.
\end{lemma}

{\bf Proof. }The Hessian matrix of $f$ is
$$\begin{pmatrix} 2x_3 & 0 & 0 & 2x_0\\
0 & 2x_3 & 0 & 2x_1 \\
0 & 0 & -2x_3 & -2x_2 \\
2x_0 & 2x_1 & -2x_2 & -6x_3
\end{pmatrix}.$$
The Hessian $H(f)$, the determinant of this matrix,
is $2^4x_3^2(x_0^2+x_1^2-x_2^2+3x_3^2)$. Therefore
$H(f)$ is negative if and only if $x_3$ is not zero and
$x_0^2+x_1^2-x_2^2+3x_3^2<0$. At these points, the Hessian matrix must have
signature $(+,-,-,-)$ or $(+,+,+,-)$. By inspecting the upper left
$2\times  2 $ submatrix of the Hessian matrix, we see that
$f$ has signature $(+,-,-,-)$ if and only if $x_3<0$
and $x_0^2+x_1^2-x_2^2+3x_3^2<0$. At such a point, we have
\begin{align*}
x_0^2+x_1^2-x_2^2-x_3^2&\leq x_0^2+x_1^2-x_2^2+3x_3^2\\
&< 0,
\end{align*}
and so $f=(x_0^2+x_1^2-x_2^2-x_3^2)x_3$ is positive. That is,
the index cone of $f$ is
$$\{(x_0,x_1,x_2,x_3)\in \R^4: x_3<0,\text{ }x_0^2+x_1^2-x_2^2+3x_3^2<0\}.$$
It is then easy to check that this index cone is nonempty (it has two connected
components, each a half of a component of a standard circular cone
in $\R^4$).

To compute the curvature of the 3-manifold $M=\{ f=1\}$ in the index cone,
it is helpful to notice that the form $f$ has a big automorphism group.
In particular, the orthogonal group of the quadratic form
$x_0^2+x_1^2-x_2^2$, fixing the coordinate $x_3$, preserves the
form $f$. This group can move any point in $\R^4$ to a point with $x_0=0$.
Since this group preserves the form $f$, it preserves the Hessian metric
associated to $f$.
Therefore it suffices to show that for every point of $M$
such that $x_0=0$, the sectional curvature is
positive on some 2-plane.

We can reduce to a lower-dimensional problem, using the symmetries of $f$.
Namely, the group $\Z/2$, acting on $\R^4$ by changing the sign
of $x_0$, preserves the form $f$, and so it preserves the metric on
the 3-manifold $\{ f=1\}$ in the index cone.
Therefore the fixed point set of this group action is
a totally geodesic submanifold. Clearly this fixed point set is
the surface $\{x_0=0, f=1\}$ in the index cone. Therefore, the sectional
curvature of this surface is equal to the sectional curvature of
the 3-manifold at the corresponding 2-plane. Thus, it suffices to show
that the surface $\{x_0=0, f=1\}$ in the index cone has positive
curvature everywhere. 

It is clear that the restriction of our metric $-1/d(d-1)\d^2 f
/\d x_i \d x_j$ on $M$ to this surface is
the analogous metric associated to the form $f|_{x_0=0}$ on $\R^3$. So we need
to show that the metric associated to the cubic form
$g:=(x_1^2-x_2^2-x_3^2)x_3$ on $\R^3$,
restricted to the surface $\{ g=1\}$, has positive curvature at all points
with $x_3<0$ and $x_1^2-x_2^2+3x_3^2<0$.

We compute that the Hessian of $g$ is $H(g)=8x_3(x_1^2-x_2^2+3x_3^2)$
and that the Clebsch invariant of $g$ is $S(g)=16$. (Recall that
this paper's definition of the Clebsch invariant makes it,
for plane cubics, equal to $2^43^4$ times
the Aronhold invariant $S$, as used in Wilson's paper \cite{Wilson}.)
By Wilson's calculation for plane
cubics (\cite{Wilson}, Theorem 5.1), generalized in Theorem \ref{curv}
in this paper,
the curvature of the metric on $g=1$ is:
\begin{align*}
K &= -\frac{9}{4}+9\frac{S(g)g^2}{H(g)^2}\\
&= -2^{-2}3^2+2^{-2}3^2\frac{(x_1^2-x_2^2-x_3^2)^2}{(x_1^2-x_2^2+3x_3^2)^2}.
\end{align*}
We are considering the region where $x_3<0$ and $x_1^2-x_2^2+3x_3^2<0$.
In this region, we have
$$x_1^2-x_2^2-x_3^2<x_1^2-x_2^2+3x_3^2<0.$$
So the above formula shows that the curvature of
the surface $g=1$ is positive everywhere
in this region. \qed

\section{Hessian metrics of constant curvature}
\label{const}

Wilson showed that for a Fermat form $f$ of any degree $d$
and any number of variables $n$, if the associated Hessian metric is
Lorentzian (that is, of signature $(1,*)$) on a nonempty open
subset of $\R^n$,
then the associated
Riemannian metric on $M=\{ f=1\}$ has constant sectional curvature
$-d^2/4$ \cite{Wilson}. In this section, we describe a larger
class of forms which give metrics of
constant curvature $-d^2/4$ on $M$. Dubrovin showed
that certain complex reflection groups (Shephard groups)
give further examples of forms with this property. In 3 variables,
Dubrovin's construction gives just one ``new'' example,
the Maschke sextic.
We ask whether
the forms we find, together with Dubrovin's example,
are the only forms in 3 variables
which give metrics of curvature $-d^2/4$ (Question
\ref{conjflat}). We reformulate this as a problem
in invariant theory (Question \ref{clebsch}), which we study for the
rest of the paper.

One justification for studying the condition that $M$ has constant
curvature $-d^2/4$ is that, at least when $n=3$ and the form
$f$ is irreducible, this is
the only possible constant value of the curvature, as one can
deduce from Theorem \ref{curv}. Another is that this condition is
equivalent to flatness of the Hessian metric on $\R^n$, by Corollary
\ref{flat}.

In contrast to the previous sections about inequalities on the curvature,
the Lorentzian condition does not play an important role here. We may as
well ask the more general question: for which real forms $f$ on $\R^n$
does the pseudo-Riemannian Hessian metric
$-1/d(d-1)\partial^2f/\partial x_i\partial x_j$ restricted to
the hypersurface $M:=\{ f=1 \}$ have constant curvature $-d^2/4$
on a nonempty open subset of $M$? As in the proof of Theorem \ref{curv},
it is convenient to work even more generally, with a holomorphic
``metric'' on a complex manifold,
meaning a nondegenerate symmetric bilinear form on the
tangent bundle. (Beware that this is not the usual kind of
metric considered in complex differential geometry, which is hermitian
rather than bilinear.) We can ask:
for which complex forms $f$ on $\C^n$
does the Hessian metric $-1/d(d-1)\partial^2f/\partial x_i\partial x_j$
restricted to the hypersurface $M:=\{ f=1\}$ have constant
curvature $-d^2/4$ on a nonempty open subset of $M$?

By Corollary \ref{flat}, the hypersurface $M$ has
constant curvature $-d^2/4$ if and only if the Hessian metric
is flat on a nonempty open subset of $\C^n$. So the above problem
is equivalent to the more natural-looking problem of classifying forms $f$
such that the Hessian metric is flat on a nonempty open subset of $\C^n$.
Since the sectional curvature is defined by algebraic formulas
(section \ref{hess}), flatness on a nonempty open subset implies flatness
on a dense open subset of $\C^n$.

Here is a simple class of flat Hessian metrics.

\begin{lemma}
\label{two}
For any form $f$ on $\C^2$, the Hessian metric associated to $f$
is flat (wherever
this makes sense, that is, where the Hessian determinant is nonzero).
\end{lemma}

{\bf Proof. }Since $M=\{f=1\}$ is a complex 1-manifold, the restriction
of the Hessian metric on $\C^2$ to $M$ has constant
curvature $-d^2/4$ on all complex 2-planes (this condition being
vacuous). By Corollary \ref{flat},
it follows that the Hessian metric on $\C^2$ is flat. \qed

In fact, this proof shows how to construct an isometry from the
Hessian metric associated to a real binary form $f$ to a standard flat
pseudo-Riemannian metric on $\R^2$. Assume, for example, that
the Hessian of $f$ has signature $(1,1)$ on some open subset of $\R^2$
where $f>0$. Then the metric we consider, $-1/d(d-1)\d ^2f/\d x_i \d x_j$,
also has signature $(1,1)$, and its restriction
to $M=\{f=1\}$ is positive definite.
Therefore we can construct an isometry from $(d^2/4)M$ to some open
subset of the
manifold $N=\{ g=1\}$, where $g(x,y):=x^2-y^2$, since $N$ is another
1-manifold with positive definite metric. By Lemma \ref{warp},
this isometry extends to an isometry from $d^2/4$ times the Hessian
metric of $f$ on $\R^2$ to the Hessian metric of $g$ on $\R^2$,
which is the standard flat Lorentzian metric on $\R^2$.

From Lemma \ref{two}, it follows that any form on $\C^3$ which can be written
$f=\alpha(x,y)+\beta(z)$ in some linear coordinate system gives
a flat pseudo-Riemannian Hessian metric, wherever the Hessian determinant
is nonzero. Indeed, the corresponding metric is the product of a flat metric
on an open subset of $\C^2$, by Lemma \ref{two}, with a metric on $\C$,
which is automatically flat.
More generally, in any dimension, this argument shows that
the Hessian metric of
any form $f=\alpha_1(x_1,x_2)+\alpha_2(x_3,x_4)+\cdots $ is flat.
This generalizes Wilson's observation (\cite{Wilson},
Introduction, Example 2) that the real Fermat form
$f=x_1^d-x_2^d-\cdots -x_n^d$ determines a metric on $M=\{ f=1\}$
with constant curvature $-d^2/4$, or equivalently a flat metric on an open
subset of $\R^n$ by Corollary \ref{flat}. In $\C^3$, we can ask if
the forms we have found, together with Dubrovin's example
of the Maschke sextic ((\cite{Dubrovinalmost},
Corollary 5.9 and Example 3),
are essentially the only ones for which the Hessian metric
is flat:

\begin{question}
\label{conjflat}
Let $f$ be a  ternary form of degree $d$ over
$\C$ whose Hessian
determinant is not identically zero. Suppose that the
Hessian metric $-1/d(d-1)\partial^2f/\partial x_i \partial x_j$
is flat on a nonempty open subset of $\C^3$. Is $f$ either
in the closure
of the set of forms which can be written as $\alpha(x,y)+\beta(z)$
in some linear coordinates, or equal to the Maschke sextic
$x^6+y^6+z^6-10(x^3y^3+y^3z^3+z^3x^3)$ in some linear coordinates?
\end{question}

By Theorem \ref{curv}, Question \ref{conjflat} is equivalent
to the following question. This question is formally
more general in that it
makes sense even for forms $f$ whose Hessian determinant is identically
zero. (Gordan and Noether showed that a ternary form whose Hessian
determinant $H(f)$ is identically zero can be written as
$\alpha(x,y)$ in some linear coordinates $(x,y,z)$ \cite{GN}.)
The question could have been asked by
the 19th-century invariant theorists. The rest
of the paper will be devoted to it.

\begin{question}
\label{clebsch}
Let $f$ be a ternary form of degree $d$ over $\C$ whose
Clebsch covariant $S(f)$ is identically zero. Is $f$ either
in the closure
of the set of forms which can be written
as $\alpha(x,y)+\beta(z)$ in some linear coordinates,
or equal to the Maschke sextic
$x^6+y^6+z^6-10(x^3y^3+y^3z^3+z^3x^3)$ in some linear coordinates?
\end{question}

\section{The closure of the set of forms $\alpha(x,y)+\beta(z)$}

To clarify Question \ref{clebsch}, in this section we give
an explicit description of the closure of the set of ternary forms
over $\C$ which can be written as $\alpha(x,y)+\beta(z)$ in some
linear coordinate system.

\begin{lemma}
\label{closure}
Let $f$ be a ternary form over $\C$. Then 
$f$ is in the closure of the set of forms which can be
written as $\alpha(x,y)+\beta(z)$ in some linear coordinates
if and only if $f$ can be written as either $\alpha(x,y)+\beta(z)$
or $\alpha(x,y)+\beta(x)z$ in some linear coordinates.
\end{lemma}

{\bf Proof. }
It is easy to show that the second condition implies the first.
That is, we have to show that any form
$f=\alpha(x,y)+\beta(x)z=\alpha(x,y)+bx^{d-1}z$ is in the closure
of the set of $GL(3,\C)$-translates of forms $\gamma(x,y)+\delta(z)$.
This is immediate by noting that
$$f=\lim_{c\arrow 0} \bigg[ \big( -\frac{x^d}{c}+\alpha(x,y)\big) +\frac{1}{c}
\big( x+\frac{cb}{d}z\big) ^d\bigg] .$$

For the converse, we have
to show that any ternary form in the closure of the set
of $GL(3,\C)$-translates of forms $\alpha(x,y)+\beta(z)$
can be written either as $\alpha(x,y)+\beta(z)$ or as $\alpha(x,y)
+\beta(x)z$, in some linear coordinates.

We use that any form $f=\alpha(x,y)+\beta(z)$ satisfies
the differential equations
$$\frac{\partial^2f}{\partial x\partial z}=\frac{\partial^2f}
{\partial y\partial z}=0.$$
Thus, in the 3-dimensional complex vector space spanned by
$\partial/\partial x$, $\partial/\partial y$, and 
$\partial/\partial z$, there is a 1-dimensional subspace
$L$ and a 2-dimensional subspace $P$ such that $LPf=0$.
Clearly this remains true for any $GL(3,\C)$-translate of $f$.
Therefore it remains true for any form $g$ in the closure
of such forms, although the line $L$ may be contained in the plane $P$.
Thus, after a change of coordinates, $g$ either satisfies
$$\frac{\partial^2g}{\partial x\partial z}=\frac{\partial^2g}
{\partial y\partial z}=0$$
or
$$\frac{\partial^2g}{\partial z^2}=\frac{\partial^2g}
{\partial y\partial z}=0.$$
In the first case we can write $g=\alpha(x,y)+\beta(z)$,
and in the second case we can write $g=\alpha(x,y)+\beta(x)z$. \qed

\section{Plane quartic curves with vanishing Clebsch covariant}

In this section, we give a positive answer to
Question \ref{clebsch} for forms
of degree at most 4.

\begin{theorem}
\label{four}
A ternary form $f$ of degree at most 4 over the complex numbers
has Clebsch covariant equal to
zero if and only if it is in the closure of the
set of forms which can be written as $\alpha(x,y)+\beta(z)$
in some linear coordinates.
\end{theorem}

{\bf Proof. }This is trivial for forms of degree at most 2 (where every form
can be written as $\alpha(x,y)+\beta(z)$). For cubic forms,
the Clebsch covariant is a constant multiple of the Aronhold
invariant $S$, and it is a standard fact from the theory of
elliptic curves that a cubic curve
has $S=0$ if and only if it is in the closure of the
set of $GL(3,\C)$-translates
of the Fermat cubic $x^3+y^3+z^3$ (\cite{DK}, Prop.\ 5.13.2).

It remains to show that a ternary quartic with Clebsch covariant zero
is in the closure of the set of $GL(3,\C)$-translates of the forms
$\alpha(x,y)+\beta(z)$. (These forms are generally not in the closure
of the orbit of the Fermat quartic $x^4+y^4+z^4$.)

Dolgachev and Kanev showed that a ternary quartic $f$ with
Clebsch covariant zero is not ``weakly nondegenerate'', in their
terminology (\cite{DK}, Cor.\ 6.6.3(iv)).
That is, in some linear coordinates,
we have either $\d^2f/\d y\d z=0$ or $\d^2 f/\d y^2=0$.
We are trying to prove a stronger conclusion with the same hypothesis,
building on their result.

Suppose we are in the first case, that is, that $\d^2 f/\d y\d z=0$.
Then we can write $f=\alpha(x,y)+\beta(x,z)$. To say more, we need
to see what the vanishing of the Clebsch covariant says about such a form.
In the notation used
to define the Clebsch covariant, the forms $b_1=f_{123}$,
$c_1=f_{223}$, and $b_2=f_{233}$ are zero. So the formula for
the Clebsch covariant becomes:
$$S(f)= (a_1a_3-a_2^2)(d_0b_0-c_0^2).$$
Thus, since the Clebsch covariant $S(f)$ is zero, we have
either $a_1a_3=a_2^2$ or $d_0b_0=c_0^2$. These two conditions 
are the same up to switching the coordinates $y$ and $z$,
so we can assume that $a_1a_3=a_2^2$, that is, $f_{113}f_{333}=f_{133}^2$.

Since we have $f=\alpha(x,y)+\beta(x,z)$, the derivative $f_3$ is 
equal to $\beta_3$, and so the above equation says that the Hessian
determinant
of $\beta_3=(\d \beta/\d z)(x,z)$ is zero. It is classical that the vanishing
of the Hessian of a binary form implies that the form is a power
of a linear form (\cite{Olver}, Prop.\ 2.23).
Thus we can write $\beta_3=(bx+cz)^{d-1}$
for some numbers $b,c$. If $c=0$, so that $\beta_3$ is a constant
multiple of $x^{d-1}$, then $\beta$ itself is a linear combination
of $x^{d-1}z$ and $x^d$. Thus, putting the $x^d$ term into $\alpha$,
we can write the form $f$ as $\alpha(x,y)+ax^{d-1}z$
for some $a$. This proves
what we want in the case $c=0$,
since the forms $\alpha(x,y)+ax^{d-1}z$
are in the closure of the set of
$GL(3,\C)$-translates of the forms $\alpha(x,y)+\beta(z)$
by Lemma \ref{closure}. It remains
to consider the case where $\beta_3=(bx+cz)^{d-1}$ with $c$ not zero.
Then $\beta$ itself is a linear combination of $(bx+cz)^d$ and $x^d$.
Putting the $x^d$ term into $\alpha$, we can write the form $f$
as $\alpha(x,y)+a(bx+cz)^d$ for some $a$. Since $c$ is not zero,
we can change coordinates to write $f=\alpha(x,y)+\beta(z)$, as we want.

It remains to consider the second case of Dolgachev and Kanev's
result, where we have $\d^2f/\d y^2=0$. Intuitively, this should be
a more special case than the previous one.
In this case, we can write
$f=\alpha(x,z)+\beta(x,z)y$. We need to work out what the vanishing
of the Clebsch covariant tells us about such a form. In the notation
defining the Clebsch covariant, the forms $c_0=f_{122}$, $c_1=f_{223}$,
and $d_0=f_{222}$ are zero. Therefore the formula for the Clebsch
covariant simplifies to:
$$S(f)=(b_0b_2-b_1^2)^2.$$
Since the Clebsch covariant $S(f)$ is zero, we have $b_0b_2=b_1^2$,
that is, $f_{211}f_{233}=f_{213}^2$. Here $f_2=\beta(x,z)$, and so
this equation says that the binary form $\beta(x,z)$ has
Hessian equal to zero. Therefore $\beta(x,z)$ is a power of a linear
form. Thus, after a linear change of coordinates in $x$ and $z$, we can write
$f=\alpha(x,z)+ax^{d-1}y$ for some $a$. This is the conclusion we want,
since such forms are in the closure of the set of $GL(3,\C)$-translates
of forms $\alpha(x,y)+\beta(z)$ by Lemma \ref{closure}. \qed

\section{An irreducible component of the set of plane curves
with vanishing Clebsch covariant}

In this section, we give further evidence for 
a positive answer to Question \ref{clebsch}.
Namely, we show that the forms which can be written
as $\alpha(x,y)+\beta(z)$ in some linear coordinates
comprise an
irreducible component of the set of all forms with vanishing
Clebsch covariant, for forms of any degree.

\begin{theorem}
\label{irred}
For any number $d$, let $Y$ be the closure of the set of ternary
forms $f$ of degree $d$ over the complex numbers which can be written as
$\alpha(x,y)+\beta(z)$ in some linear coordinates. Then $Y$
is an irreducible component of the set of forms $f$ with
Clebsch covariant $S(f)=0$. Moreover, the scheme defined by
$S(f)=0$ is reduced at a general point of $Y$.
\end{theorem}

{\bf Proof. }For $d$ at most 2, every ternary form of degree $d$
belongs to $Y$, and the Clebsch covariant is identically zero;
so the statement is clear. For $d=3$, 
we know the stronger statement (Theorem \ref{four}) that
the equation $S(f)=0$ defines the set $Y$. Also, $S(f)$ is an
irreducible invariant on the space of plane cubic curves, and so it
defines $Y$ as a reduced hypersurface.
So we can assume
that $d$ is at least 4.

It is clear that the set $Y$ is irreducible.
We will show that there is a point of $Y$ in a neighborhood of which
the equation
$S(f)=0$ defines $Y$ as a smooth scheme. This implies the same statement
on a dense open subset of $Y$, and hence
implies both statements of the theorem.

Clearly, a general point $f$
of $Y$ can be written, in some linear coordinates, as
$f=\alpha(x,y)+\beta(z)$,
for some forms $\alpha(x,y)$ and $\beta(z)$ of degree $d$. Thus we need
to choose forms $\alpha(x,y)$ and $\beta(z)$ such that
the equation $S(f)=0$ defines $Y$ as a smooth scheme in a neighborhood
of the point $f$. We will take $\beta(z)=z^d$. (This is no loss
of generality, since a general form $\beta(z)=cz^d$ can be put
into this form by scaling $z$.)

There is an easy lower bound for the Zariski tangent space of $Y$
at the point $f$. First, we can vary the forms $\alpha(x,y)$ and
$\beta(z)$; and then, we can move $f=\alpha(x,y)+\beta(z)$ by
elements of $GL(3,\C)$, in the direction of an element
of the Lie algebra $\gl(3,\C)$. The Lie subalgebra $\gl(2)\times \gl(1)$
maps the vector space of forms $\gamma(x,y)+\delta(z)$ into itself,
and so it suffices to consider the 4-dimensional vector space
$\gl(3)/\gl(2)\times \gl(1)$. It is spanned by the four infinitesimal
transformations $x\mapsto x+\epsilon z$, $y\mapsto y$, $z\mapsto z$;
and likewise $y\mapsto y+\epsilon z$ with other variables unchanged;
and likewise $z\mapsto z+\epsilon x$; and likewise $z\mapsto z+\epsilon y$.
These transformations move the form $f=\alpha(x,y)+z^d$ in the
directions of the following 4 forms: $\alpha_1z$, $\alpha_2z$,
$xz^{d-1}$, $yz^{d-1}$. (Following section \ref{conventions},
$\alpha_1$ means $\d \alpha/\d x$ and $\alpha_2$ means $\d \alpha
/\d y$ here.) Thus the Zariski tangent space to $Y$
at the point $f$ is at least the span of the forms $x^d,x^{d-1}y,\ldots,y^d$,
$\alpha_1z$, $\alpha_2z$,
$xz^{d-1}$, $yz^{d-1}$, $z^d$. Clearly, for a general form $\alpha(x,y)$
of degree $d$, the forms listed are linearly independent.

The theorem will be proved if we can show that there is a form
$\alpha(x,y)$ of degree $d$ such that
the Zariski tangent space to the scheme
$S(f)=0$ at the point $f=\alpha(x,y)+z^d$ is spanned by the above forms.
To compute the Zariski tangent space to the scheme $S(f)=0$,
we need to compute $S(f+\epsilon g) \pmod{\epsilon^2}$
for an arbitrary form $g(x,y,z)$, where $f=\alpha(x,y)+z^d$.
From the definition of the Clebsch covariant, we compute:
\begin{multline*}
S(\alpha(x,y)+z^d+\epsilon g)=\epsilon d^4(d-1)^4(d-2)^4
[g_{113}(\alpha_{112}\alpha_{222}-\alpha_{122}^2)\\
+g_{123}(\alpha_{112}\alpha_{122}-\alpha_{111}\alpha_{222})+g_{223}
(\alpha_{111}\alpha_{122}-\alpha_{112}^2)] \pmod{\epsilon^2}.
\end{multline*}
Thus, we need to show that there is a binary form $\alpha(x,y)$
of degree $d$ such that the vanishing of the expression in brackets for
a ternary form $g(x,y,z)$ of degree $d$ implies that $g$
is in the span of the forms $x^d,x^{d-1}y,\ldots,y^d$,
$\alpha_1z$, $\alpha_2z$,
$xz^{d-1}$, $yz^{d-1}$, $z^d$.

It is convenient to observe that the expression in brackets sends
the part of the polynomial $g$ such that $z$ has a given exponent $a$
to a polynomial such that $z$ has exponent $a-1$. As a result,
we can consider the problem separately for each part of $g$
of the form $h(x,y)z^a$. Thus the question becomes one about
binary forms only. Namely, for binary forms $\alpha(x,y)$ and $h(x,y)$,
define
$$T(\alpha,h)=h_{11}(\alpha_{112}\alpha_{222}-\alpha_{122}^2)+h_{12}
(\alpha_{112}\alpha_{122}-\alpha_{111}\alpha_{222})+h_{22}
(\alpha_{111}\alpha_{122}-\alpha_{112}^2).$$
This is an $SL(2)$-equivariant differential operator in $\alpha$ and
$h$. We need to show that for any $d\geq 4$, there is a form
$\alpha(x,y)$ of degree $d$ with the following properties.
First, for any form $h(x,y)$ of degree $r$ with
$2\leq r\leq d-2$, if $T(\alpha,h)=0$ then $h=0$. Second,
any form $h(x,y)$ of degree $d-1$ with $T(\alpha,h)=0$ must be a linear
combination of the derivatives $\alpha_1$ and $\alpha_2$.

In fact, we can prove even more, as follows. This will complete
the proof of Theorem \ref{irred}. It would be preferable
to have a more geometric interpretation of the operator $T(\alpha,h)$,
but it happens that we can get by without that.

\begin{lemma}
\label{general}
Let $\alpha(x,y)$ be a very general binary form of degree $d\geq 4$
(that is, a form outside countably many proper subvarieties
of the space of all forms of degree $d$).
Then the binary forms $h(x,y)$ of any degree such that $T(\alpha,h)=0$
are the linear combinations of: $1$ in degree 0, $x$ and $y$ in degree 1,
the derivatives $\alpha_1$ and $\alpha_2$ in degree $d-1$,
and $\alpha$ in degree $d$.
\end{lemma}

{\bf Proof. }First, it is a straightforward calculation that
for any form $\alpha(x,y)$, the operator $T(\alpha,h)$ vanishes
when $h$ has degree at most 1, and also that
$T(\alpha,\alpha_1)=T(\alpha,\alpha_2)=T(\alpha,\alpha)=0$. The operator
$T(\alpha,h)$ is linear in $h$, and so $T(\alpha,h)$ vanishes when
$h$ is any linear combination of the forms $1,x,y,\alpha_1,\alpha_2,\alpha$.

To prove the lemma for a very general form $\alpha(x,y)$,
it suffices to prove it for a single form $\alpha(x,y)$ of degree $d$.
(We get the conclusion only for $\alpha$ very general, that is,
outside countably many proper subvarieties, because the statement
involves forms $h$ of arbitrary degree. For the application to
Theorem \ref{irred}, we only need Lemma \ref{general}
for forms $h$ of degree at most $d-1$, and that weakened form
of the Lemma holds for forms $\alpha$ outside only finitely many
proper subvarieties.)

We take $\alpha(x,y)=\binom{d}{2}x^{d-2}y^2$. Then, for any form $h$,
we compute that
$$T(\alpha,h)=d^2(d-1)^2(d-2)^2x^{2d-8}[x^2h_{11}-(d-3)xyh_{12}
+\frac{1}{2}(d-2)(d-3)y^2h_{22}].$$
Clearly a form $h$ has $T(\alpha,h)$ equal to zero
if and only if the expression
$U(\alpha,h)$ in brackets is zero. We compute that
$$U(\alpha,x^iy^j)=x^iy^j(i(i-1)-(d-3)ij+\frac{1}{2}(d-2)(d-3)j(j-1)).$$
Thus the differential operator $h\mapsto U(\alpha,h)$ is diagonalized
on the basis of monomials $x^iy^j$. It follows that the vector space
of forms $h$ with $T(\alpha,h)$ equal to zero is spanned by the
set of monomials $x^iy^j$ such that
$$i(i-1)-(d-3)ij+\frac{1}{2}(d-2)(d-3)j(j-1)=0.$$

For fixed $j$, let us view this equation as a quadratic equation
for $i$. As such, its discriminant $b^2-4ac$ is
$$\Delta=1-(d-1)(d-3)j(j-2).$$
Since $d$ is at least 4, $(d-1)(d-3)$ is at least 3. So we read off
that $\Delta$ is negative for $j$ at least 3. Thus, there are no
integral (or even real) solutions $i$ unless $j\leq 2$. Solving our
quadratic equation for $j=0,1,2$ gives that the only solutions
$(i,j)$ in natural numbers are: $(0,0)$, $(1,0)$, $(0,1)$, $(d-2,1)$,
$(d-3,2)$, and $(d-2,2)$. Thus, for $\alpha=\binom{d}{2}x^{d-2}y^2$
as we have been considering, the vector space of forms $h$ such that
$T(\alpha,h)=0$ is spanned by $1,x,y,x^{d-2}y, x^{d-3}y^2$,
and $x^{d-2}y^2$; equivalently, it is spanned by $1,x,y,
\alpha_1,\alpha_2$, and $\alpha$. This proves the lemma for the
particular form $\alpha=\binom{d}{2}x^{d-2}y^2$. As we have said,
this implies the lemma for very general forms $\alpha$. \qed


\small \sc DPMMS, University of Cambridge, Wilberforce Road,
Cambridge CB3 0WB, England

b.totaro@dpmms.cam.ac.uk
\end{document}